\documentclass{elsarticle}
  
\usepackage[utf8]{inputenc}
  
\usepackage{amsmath}
\usepackage{bm}
\usepackage{tikz}
\usepgflibrary{arrows}
\usepackage{rotating}  
\usepackage{comment,xspace}
\graphicspath{{pictures/}}
 
\usepackage{hyperref}

\newproof{pf}{Proof}

\journal{Journal of Computational and Applied Mathematics  }

\begin{document}

\begin{frontmatter}

\title{Computational identification of adsorption and desorption parameters for pore scale transport in periodic porous media}

\author[nefu]{Vasiliy V. Grigoriev}
\ead{d$\_$alighieri@rambler.ru}

\author[itwm,ban]{Oleg~Iliev\corref{cor}}
\ead{iliev@itwm.fraunhofer.de}

\author[nsi,rudn]{Petr N. Vabishchevich}
\ead{vabishchevich@gmail.com}

\address[nefu]{Multiscale Model Reduction Laboratory, North-Eastern Federal University, 58 Belinskogo st, Yakutsk, 677000, Russia}
\address[itwm]{Fraunhofer Institute for Industrial Mathematics ITWM,  Fraunhofer-Platz 1, D-67663 Kaiserslautern, Germany}
\address[ban]{Institute of Mathematics and Informatics, Bulgarian Academy of Science, Sofia, Bulgaria}
\address[nsi]{Nuclear Safety Institute, Russian Academy of Sciences, 52, B. Tulskaya, Moscow, Russia}
\address[rudn]{RUDN University, 6 Miklukho-Maklaya st, Moscow, 117198, Russia}

\cortext[cor]{Corresponding author}

\begin{abstract}
Computational identification of unknown adsorption and desorption rates is discussed in conjunction with reactive flow considered at pore scale. The reactive transport is governed by incompressible Stokes equations, coupled with convection-diffusion equation for species transport. The surface reactions, namely adsorption and desorption, are accounted via Robin boundary conditions. Henry and Langmuir isotherms are considered. Measured concentration of the specie at the outlet of the domain has to be provided to carry out the identification procedure. Deterministic and stochastic parameter identification approaches are considered. The influence of the noise in the measurements on the accuracy of the identified parameters is discussed. Multistage identification procedure is suggested for the considered class of problems. The proposed identification approach is applicable for different geometries (random and periodic) and for a range of process parameters. In this paper the potential of the approach is demonstrated in identifying parameters of Langmuir isotherm for low Peclet and low Damkoler numbers reactive flow in a 2D periodic porous media with circular inclusions. Simulation results for random porous media and other regime parameters are subject of follow up papers. Finite element approximation in space and implicit time discretization are exploited. 
\end{abstract}

 
\end{frontmatter}

\section{Introduction}


The reactive transport in porous media is important component of many industrial and environmental problems like water purification, soil pollution and remediation, catalytic filters, $CO_2$ storage, oil recovery, etc., to name just a few. Historically, most of the theoretical and experimental research on transport in porous media in general, and on reactive transport in particular, has been carried out at macroscopic Darcy scale \cite{bear2013dynamics,helmig1997multiphase}. In many cases, the bottleneck in performing computational modeling of reactive transport is the absence of data for the pore scale adsorption and desorption rate (or in more general, the parameters of the heterogeneous reactions). Despite the progress in developing devices to perform experimental measurements at the pore–scale, experimental characterization of these rates is still a very challenging task. 

In the case of heterogeneous (surface) reactions at pore scale, the species transport is coupled to surface reaction via boundary conditions. When the reaction rates are not known, their identification falls into the class of boundary value inverse problems, \cite{lavrentev1986ill,alifanov2011inverse,isakovinverse,samarskii2007numerical}. The additional information which is needed to identify the parameters is often provided in form of dynamic change of the  concentration at the outlet (e.g., so called breakthrough curves).  In the literature, inverse problems for porous media flow are discussed mainly in connection with parameter identification for macroscopic, Darcy scale problems. An overview on inverse problems in groundwater Darcy scale modeling can be found in  \cite{sun2013inverse}. Identification of parameters for pore scale models is discussed in this paper, and the algorithms from \cite{sun2013inverse} and other papers discussing parameter identification at macroscale can not be applied here without modification. Let us shortly mention some general approaches for solving inverse problems.  

Different algorithms can be applied for solving parameter identification problems, see, e.g. 
\cite{tarantola2005inverse,Aster2013}. Many of the algorithms exploit deterministic methods based on Tikhonov regularization technique  \cite{tikhonov1977solutions,engl2014inverse} and target at minimizing a functional of the difference between measured and computed quantities. An important part of such algorithms is the definition of feasible set of parameters on which the functional is minimized. Local or global optimization procedures are used in the optimization \cite{horst2013handbook,nocedal2006numerical}.
In this sense it could be pointed out that there is certain similarity between the mathematical formulation of an optimization problem and of a parameter identification one.

Stochastic-deterministic methods are also popular approach for solving parameter identification problems. A variant of the method based on deterministic sampling of points looks appropriate for the topic considered here. A stochastic approach for global optimization in its simplest form consists only of a random search and it is called Pure Random Search \cite{Zhigljavsky2008}. In this case the residual functional is evaluated at randomly chosen points from the feasible set. 
Sobol sequences \cite{sobol1976uniformly,sobol1979systematic} can be used for sampling. The sensitivity analysis tool SALib \cite{Herman2017} has shown to be appropriate tool for this. Such an approach is successfully used, for example, in multicriteria parameter identification \cite{sobol1981choosing}. 

The solution of the inverse problems we are interested in, is composed of two ingredients: (multiple) solution of the direct (called also forward) problem, and the parameter identification algorithm.

Note that the considered here parameter identification approach is different from approaches like determining the single fibre (or collector) efficiency. Originally \cite{Baron} this additive approach used analytical solution for flow and transport around a single sphere or cylinder to evaluate the contribution of the convection and diffusion to the particles deposition and to propose a formula for this deposition depending on Pe number and fiber diameter. Later on several modification were proposed, for example in \cite{Mesina} a convection diffusion equation with zero Dirichlet boundary condition on the surface of a single fiber, in order to calculate the specie flux to the surface and to fit an algebraic formula for the single fiber efficiency. In \cite{Wood} it was stated that zero Dirichlet boundary condition on the surface of the pore is simplification and introduced a jump described by Henry isotherm. The adsorption and desorption parameters in this isotherm were derived from the solution of Smoluchovski equation at nanoscale. This was further embedded in volume averaging upscaling procedure from nanoscale to Darcy scale. As mentioned above, our aim is to have an identification procedure applicable to any periodic and random porous media, for various process regimes, various kinetics of the heterogeneous reactions, and thus to go beyond the limited applicability of formula like single collector efficiency. On the other hand, we assume that it is known which isotherm has to be used for the processes we model, and we do not plan to go to nanoscale to derive isotherms.  

The goal of this paper is to contribute to the understanding of the formulation and the solution of a class of parameter identification problems for pore scale reactive transport in the case when the measured concentration of the specie at the outlet of the domain is provided as extra information in order to carry out the identification procedure. Deterministic and stochastic parameter identification approaches are considered. The influence of the noise in the measurements on the accuracy of the identified parameters is discussed. Multistage identification procedure is suggested for the considered class of problems. The proposed identification approach is applicable for different geometries (random and periodic) and for a range of process parameters. In this paper the potential of the approach is demonstrated in identifying parameters of Langmuir isotherm for low Peclet and low Damkoler numbers reactive flow in a 2D periodic porous media with circular inclusions. It is supposed that this paper is the first one in series of papers dedicated to this topic. Simulation results for random porous media and other regime parameters are subject of follow up papers. 

The reminder of the paper is organized as follows.
The direct problem is considered in Section 2. At pore scale, single phase laminar flow described by incompressible Stokes equations, and solute transport described by convection-diffusion equation, are considered. The surface reaction is accounted for in the boundary conditions. Henry and Langmuir adsorption isotherms are considered here \cite{kralchevsky2008chemical}, identification is carried out for  adsorption and desorption parameters in the Langmuir isotherm. Section 3 is dedicated to description of the used computational algorithm. Finite element method is exploited after triangulation of the computational domain. The numerical investigation of the grid convergence and sensitivity with respect to parameters is also presented in this Section. The set up of the parameter identification problem for reactive flow in porous media is described in Section 4. Computational algorithm for statistical evaluation of adsorption parameters, based on multiple solution of direct problems, is presented there. 
Finally, Section 5 summarizes the results presented in this paper.

\section{Mathematical model} 

\begin{figure}[ht] 
  \begin{center}
    \begin{tikzpicture}
       \shade[top color=blue!5, bottom color=blue!5] (7,0) rectangle +(3,5);
       \shade[left color=blue!50, right color=blue!5] (3,0) rectangle +(4,5);
       \shade[top color=blue!50, bottom color=blue!50] (0,0) rectangle +(3,5);
       \draw [dashed] (0, 0) rectangle +(10,5);
       \draw [color=red] (0,4.25) rectangle +(10,0.5);
       \foreach \y in {3,...,6} 
	\foreach \x in {0,...,4} {
  	\filldraw [fill=green,draw=black] (\y,\x+0.25) circle (0.24);
  	\filldraw [fill=green,draw=black] (\y+0.5,\x+0.75) circle (0.24);
       }
       \draw [<-, line width=2, color=blue] (1,1.5) -- (-1,1.5);
       \draw [<-, line width=2, color=blue] (1,2.5) -- (-1,2.5);
       \draw [<-, line width=2, color=blue] (1,3.5) -- (-1,3.5);
       \draw [<-, line width=2, color=blue] (11,1.5) -- (9,1.5);
       \draw [<-, line width=2, color=blue] (11,2.5) -- (9,2.5);
       \draw [<-, line width=2, color=blue] (11,3.5) -- (9,3.5);
       \draw [-] (2,4.5) -- (2,6);
       \draw  (2,6.4) node {$\Omega_f$};
       \draw [-] (4,4.375) -- (4,6);
       \draw  (4,6.4) node {$\Omega_s$};  
       \draw [-] (6,4.5) -- (6,6);
       \draw  (6,6.4) node {$\Gamma_s$}; 
       \draw [-] (0,4.5) -- (-0.5,6);
       \draw  (-0.5,6.4) node {$\Gamma_{in}$};
       \draw [-] (7.5,4.25) -- (8,6);
       \draw [-] (8.5,4.75) -- (8,6);
       \draw  (8,6.4) node {$\Gamma_{sim}$};
       \draw [-] (10,4.5) -- (10.5,6);
       \draw  (10.5,6.4) node {$\Gamma_{out}$};  
       \draw [->] (-0.5,-0.5) -- (1.5,-0.5);
       \draw  (1.9,-0.5) node {$x_1$};  
       \draw [->] (-0.5,-0.5) -- (-0.5,0.5);
       \draw  (-0.5,0.9) node {$x_2$};  
    \end{tikzpicture}
    \caption{Sketch of the pore scale domain} 
   \label{f-1}
  \end{center}
\end{figure}
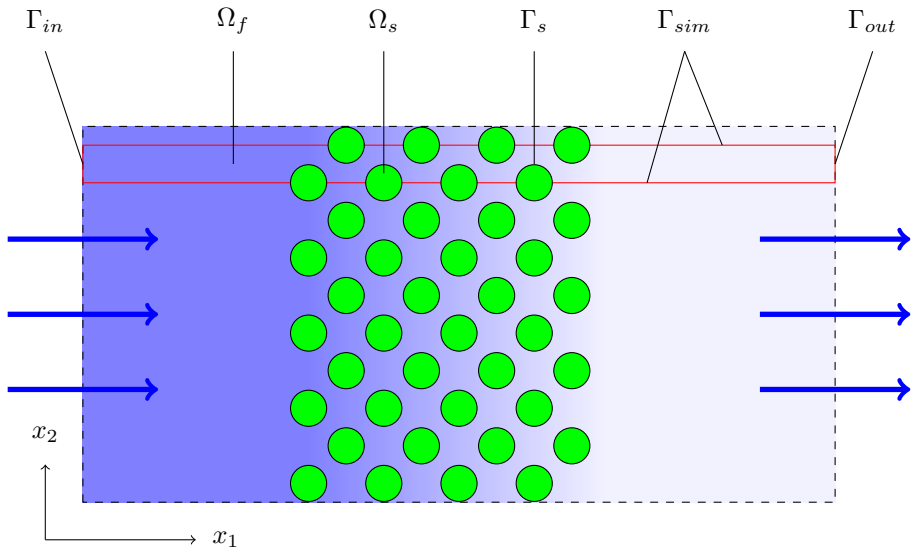

Two dimensional pore scale transport of a dissolved substance in the presence of adsorption and desorption is considered. The coordinate along the porous domain (the vertical one) is denoted as $x_2$, while the coordinate across the porous domain (the horizontal one, coinciding with mean velocity direction) is denoted as $x_1$). Sketch of the flow is shown on Fig.\ref{f-1}.
The porous media geometry is modeled as periodic arrangement of cylinders.
A part of the domain, namely $\Omega_f$ is occupied by a fluid, while the other part is occupied be the obstacles $\Omega_s$. The obstacle surfaces (where the reaction occurs) are denoted by $\Gamma_s$, while symmetry lines are denoted by $\Gamma_{sim}$. It is supposed that dissolved substance is introduced via the inlet boundary 
$\Gamma_{in}$, and the part of the substance which did not react flows out via $\Gamma_{out}$. The computational domain consists from $\Omega_f$. 

\subsection{Flow problem} 

The flow in the pores is often slow, it is described here by the steady state incompressible Stokes equations:
\begin{equation}\label{1}
\nabla p - \mu \nabla^2 \bm{u} = 0,
\end{equation} 
\begin{equation}\label{2}
 \nabla \cdot \bm{u} = 0, 
 \quad  \bm{x} \in \Omega_f, 
 \quad  t > 0,
\end{equation} 
where $\bm{u}(\bm{x})$ and $p(\bm{x})$ are the fluid velocity and pressure, respectively,
while $\mu > 0$ and $\rho > 0$ are the viscosity and the density, which we assume to be constants \cite{bear2013dynamics,acheson2005elementary}. 


Denote by $\bm n$ the outer normal vector to the boundary. 
Suitable boundary conditions on $\partial \Omega_f$ are specified. 
The velocity of the fluid $\bar{u}$ is prescribed at the inlet 
\begin{equation}\label{3}
 \bm{u} \cdot \bm n = \bar{u},
 \quad \bm{u} \times  \bm n = 0,
 \quad \bm x \in \Gamma_{in} .
\end{equation} 
At the outlet, pressure and absence of tangential force are prescribed 
\begin{equation}\label{4}
 p - \bm \sigma \bm n \cdot \bm n = \bar{p},
 \quad \bm \sigma \bm n \times \bm n = 0, 
 \quad \bm x \in \Gamma_{out} . 
\end{equation}  
Standard no-slip and no-penetration conditions are prescribed on the solid walls:
\begin{equation}\label{5}
 \bm{u} \cdot \bm n = 0,
 \quad \bm{u} \times  \bm n = 0,
 \quad \bm x \in \Gamma_{s} .
\end{equation}
Symmetry conditions are prescribed on the symmetry boundary of the computational domain:
\begin{equation}\label{6}
 \bm{u} \cdot \bm n = 0,
 \quad \bm \sigma \bm n \times \bm n = 0, 
 \quad \bm x \in \Gamma_{sim} .  
\end{equation} 

\subsection{Species Transport} 

The concentration of the solute in the fluid is denoted by $c(\bm x, t)$.
The unsteady solute transport in absence of homogeneous reactions is governed by convection diffusion equation
\begin{equation}\label{7}
 \frac{\partial c }{\partial t} + \nabla (\bm{u} c) 
 - D \nabla^2 c = 0 ,
 \quad  \bm{x} \in \Omega_f, 
 \quad  t > 0,
\end{equation} 
where $D > 0$ is the solute diffusion coefficient which is assumed to be scalar and constant. 

The concentration of the solute at the inlet is assumed to be known:
\begin{equation}\label{8}
 c(\bm x, t) = \bar{c},
 \quad \bm x \in \Gamma_{in} ,
\end{equation} 
where $\bar{c} > 0$ is assumed to be constant.
Zero diffusive flux of the solute at the outlet and on
the external boundaries of the domain is prescribed as follows:
\begin{equation}\label{9}
 D  \nabla c \cdot \bm{n} = 0,
 \quad \bm x \in \Gamma_{sim} \cup  \Gamma_{out} .
\end{equation} 
Note that convective flux via the outlet is implicitly allowed by the above equations.
The surface reactions that occur at the obstacles' surface $\Gamma_s$ satisfy the mass conservation law, in this particular case meaning that the change in adsorbed surface concentration is equal to the flux from the fluid to the surface. This is described as
\begin{equation}\label{10}
 \frac{\partial m}{\partial t} = - D  \nabla c \cdot \bm{n} ,
 \quad \bm x \in \Gamma_{s},
\end{equation} 
where $m$ is the surface concentration of adsorbed solute \cite{kralchevsky2008chemical}.
A mixed kinetic–diffusion adsorption description is used:
\begin{equation}\label{11}
 \frac{\partial m}{\partial t} = f(c, m) .
\end{equation} 
For reactive boundaries, the choice of $f$  and its dependence on $c$ and $m$ is critical for a correct description of the reaction dynamics 
at the solid–fluid interface. A number of different isotherms (i.e., different functions $f(c,m)$) exist for describing these dynamics, 
dependent on the solute attributes, the order of the reaction, and the interface type.

The simplest of these is the Henry isotherm, which assumes a linear relationship between the near surface
concentration and the surface concentration of the adsorbed particles, and takes the form
\begin{equation}\label{12}
 f(c, m) = k_a c - k_d m,
\end{equation} 
Here $k_a \geq 0$ is the rate of adsorption, measured in unit length per unit time, 
and $k_d \geq 0$ is the rate of desorption, measured per unit time.
The Langmuir adsorption isotherm is a more complicated, three parameter model:
\begin{equation}\label{13}
 f(c, m) = k_a c \left (1 - \frac{m}{m_\infty} \right )  - k_d m .
\end{equation} 
Here $m_\infty > 0$  is the maximal possible adsorbed surface concentration. 
In comparison to the Henry isotherm (\ref{12}), the
Langmuir isotherm (\ref{13}) predicts a decrease in the rate of adsorption as the adsorbed concentration increases
due to the reduction in available adsorption surface. 

The formulation of the initial -- boundary value problem in addition to the governing equations  
(\ref{1}), (\ref{2}), (\ref{7}), boundary conditions (\ref{3})--(\ref{6}), (\ref{8})--(\ref{11}), and specified isotherm (\ref{12}) or (\ref{13}), requires specification of initial conditions:
\[
 \bm{u}(\bm x, 0) = \bm{\bar{u}}(\bm x), 
 \quad \bm x \in \Omega_f , 
\]
\begin{equation}\label{14}
 c(\bm x, 0) = c_0(\bm x), 
 \quad \bm x \in \Omega_f . 
\end{equation} 
\begin{equation}\label{15}
 m(\bm x, 0) = m_0(\bm x), 
 \quad \bm x \in \Gamma_s. 
\end{equation} 


\subsection{Dimensionless form of the equations} 

When a problem like the above one need to be solved for a range of parameters (what is our goal here),
working with dimensionless form of the equations give definitive advantages.
For the dimensionless variables (velocity, pressure, concentration) below, the same notations are used as for the dimensional ones. The height of the computational domain $\Omega_f$, namely $l$, is used for scaling spatial sizes, the scaling of the velocity is done by the inlet velocity $\bar{u}$, and the scaling of the concentration is done by the inlet concentration $\bar{c}$.

The Stokes Eq.(\ref{3}) and its boundary conditions in dimensionless remain unchanged, keeping mind that in this case they are written with respect to dimensionless velocity and pressure, and considering the dimensionless viscosity to be equal to one.

In dimensionless form Eq. (\ref{7}) reads
\begin{equation}\label{20}
 \frac{\partial c }{\partial t} + \nabla (\bm{u} c) 
 - \frac{1}{\mathrm{Pe}}  \nabla^2 c = 0 ,
 \quad  \bm{x} \in \Omega_f, 
 \quad  t > 0,
\end{equation} 
where
\[
 \mathrm{Pe} =  \frac{l \bar{u} }{D} 
\] 
is the Peclet number.

Further on, Eq. (\ref{8}) is transformed into 
\begin{equation}\label{21}
 c = 1,
 \quad \bm x \in \Gamma_{in} ,
\end{equation} 
while the boundary condition (\ref{9}) take the form
\begin{equation}\label{22}
 \nabla c \cdot \bm{n} = 0,
 \quad \bm x \in \Gamma_{sim} \cup  \Gamma_{out} .
\end{equation} 
The dimensionless form of Eq.(\ref{10}) is given by
\begin{equation}\label{23}
 \frac{\partial m}{\partial t} = - \nabla c \cdot \bm{n} ,
 \quad \bm x \in \Gamma_{s},
\end{equation}
where $m$ is scaled as follows:
\[
 \bar{m} =  l \bar{c} .
\]  
In dimensionless form the adsorption relations, in the case of Henry isotherm, are written as follows
\begin{equation}\label{24}
 \frac{\partial m}{\partial t} = \mathrm{Da}_a c - \mathrm{Da}_d m,
 \quad \bm x \in \Gamma_{s},
\end{equation}
where the adsorption and desorption Damkoler numbers are given by
\[
 \mathrm{Da}_a = \frac{k_a}{\bar{u}} ,
 \quad  \mathrm{Da}_d = \frac{k_d l}{\bar{u}} .
\] 
In the case when we consider Langmuir isotherm, (\ref{11}), (\ref{13}), the following dimensionless relation is used
\begin{equation}\label{25}
 \frac{\partial m}{\partial t} = \mathrm{Da}_a c \left (1- \frac{m}{\mathrm{M}} \right )   - \mathrm{Da}_d m,
 \quad \bm x \in \Gamma_{s},
\end{equation}
where the dimensionless parameter $M$ is given by 
\[
 \mathrm{M} = \frac{m_\infty }{l \bar{c}} .
\] 


\section{Numerical solution of the direct problem} 

Finite Element Method, FEM, is used for space discretization of the above problem, together with implicit discretization in time.
The algorithm used here for solving the direct problem is practically identical of the algorithm used to study oxidation in \cite{oxidation}. Therefore some details, e.g., grid convergence studies, will be omitted here.

\subsection{Geometry and grid} 

The computational domain is a rectangle with a dimensionless height of $x_2=1$ and dimensionless length of $x_1=17.5$, in which ten half cylinders are embedded. 
The distance between the centers of the cylinders in $x_1$ direction is 1.5 dimensionless units, the radius of cylinders is 0.4 dimensionless units.

The computational domain $\Omega_f$ is triangulated using the grid generator Gmsh (website gmsh.info) \cite{Gmsh}. 
The script for preparing the geometry is written in Python. In order to control the accuracy of results, computations on consecutively refined grids are performed.

\subsection{Computation of steady state single phase fluid flow}

One way coupling is considered here. The fluid flow influences the species transport, but there is no back influence of the species concentration on the fluid flow. Based on this, the flow is computed in advance. The FEM approximation of the steady state flow problem \cite{gresho200incompressible} is based on variational formulation of the considered boundary value problem (\ref{1}), (\ref{2}), (\ref{3}), (\ref{4}), (\ref{5}), (\ref{6}). 
The following functional space $\bm V$ is defined for the velocity $\bm u$  ($\bm u \in \bm V$):
\[
\begin{split}
 \bm V = \{ \bm v \in \bm H^1(\Omega) : \  & \bm{u} \cdot \bm n = 1, \
 \bm{u} \times  \bm n = 0 \ \mathrm{on} \ \Gamma_{in} , 
 \\ & \bm{u} = 0 \ \mathrm{on} \ \Gamma_{s}, \ \bm{u} \cdot \bm n = 0 \
 \ \mathrm{on} \ \Gamma_{sim} \} .
\end{split} 
\] 
Test function $\bm v \in \hat{\bm V}$, where
\[
 \hat{\bm V} = \{ \bm v \in \bm H^1(\Omega_f) : \ \bm{u} = 0 \ \mathrm{on} \ \Gamma_{in} , \ 
 \bm{u} = 0 \ \mathrm{on} \ \Gamma_{s}, \ \bm{u} \cdot \bm n = 0 \
 \ \mathrm{on} \ \Gamma_{sim} \} .
\] 
For the pressure $p$ and the related test functions $q$, it is required that $p, q \in Q$, where 
\[
 Q = \{ g \in L_2(\Omega_f) : \ q = 0 \ \mathrm{on} \  \Gamma_{out} \} .
\] 

Let us multiply Eq.(\ref{1}) by $\bm v$, Eq.(\ref{2}) by $q$, and integrate over the computational domain. Taking into account the boundary conditions (\ref{3}), (\ref{4}), (\ref{5}), (\ref{6}), the following system of equations is obtained
with respect to $\bm v \in \bm V$, $q \in Q$
\begin{equation}\label{26}
 a(\bm u, \bm v) - b(\bm v, p) = 0 \ \forall \bm v \in \hat{\bm V} ,
\end{equation} 
\begin{equation}\label{27}
 b(\bm u, q) = 0 \ \forall q \in Q .
\end{equation} 
Here
\[
  a(\bm u, \bm v) := \int_{\Omega_f} \nabla \bm u \cdot \nabla \bm v \, d \bm x , 
\] 
\[
 b(\bm v, p) := \int_{\Omega_f}(\nabla \cdot \bm v) q \, d \bm x . 
\] 

For the FEM approximation of the velocity, the pressure, and the respective test functions, the following finite dimensional subspaces are selected
$\bm V_h \subset \bm V$, $\hat{\bm V}_h \subset \hat{\bm V}$
and $Q_h \subset Q$.
Taylor-Hood $P_2-P_1$ elements \cite{taylor1973numerical} are used here.
These are continuous $P_2$ Lagrange elements for the velocity components 
and continuous $P_1$ Lagrange elements for the pressure field.
The computations are carried out using the computing platform for partial differential equations FEniCS (website fenicsproject.org) \cite{LoggMardalEtAl2012a,AlnaesBlechta2015a}.

\begin{sidewaysfigure}
  \begin{center}
    \includegraphics[width=0.5\linewidth] {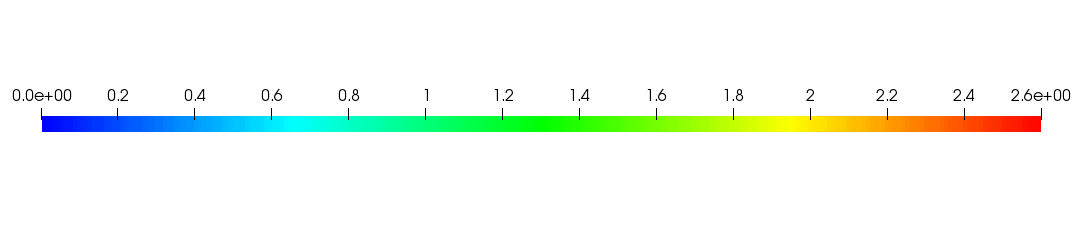} \vspace{-8mm} \\
    \includegraphics[width=1\linewidth] {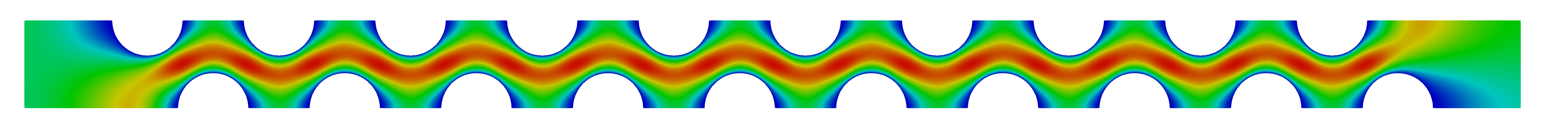} \\
	a \\
    \includegraphics[width=0.5\linewidth] {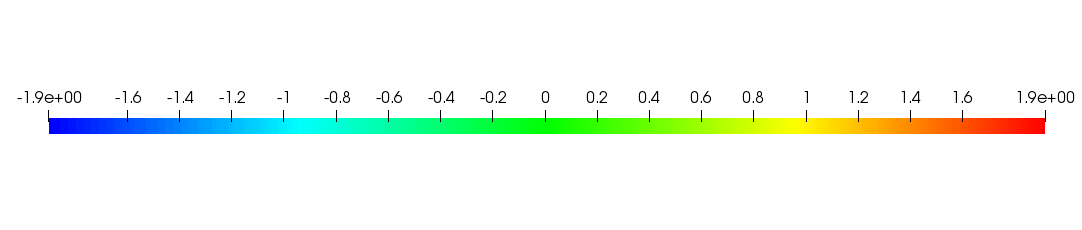}  \vspace{-8mm} \\
    \includegraphics[width=1\linewidth] {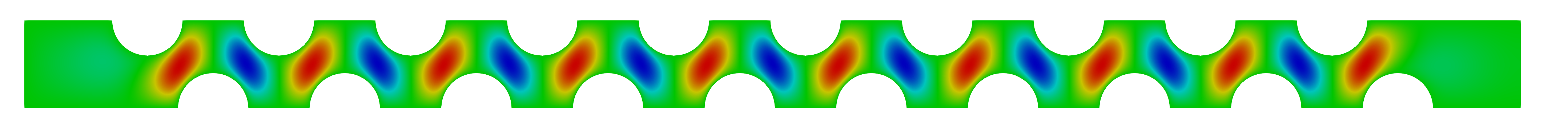} \\
	b \\
    \includegraphics[width=0.5\linewidth] {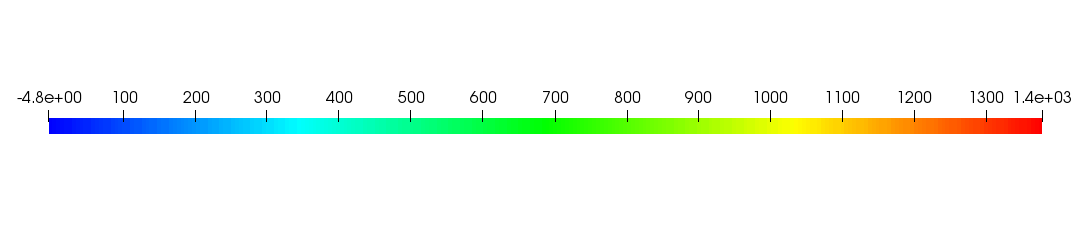}  \vspace{-8mm} \\
    \includegraphics[width=1\linewidth] {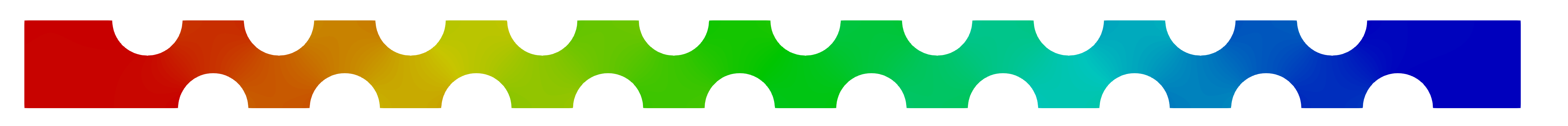} \\
	c 
  \caption{Velocity and pressure for $\mathrm{Re} = 1$: a --- horizontal component of velocity vector $u_1$,
           b --- vertical component of velocity vector $u_2$, c ---pressure}
  \label{f-2}
  \end{center}
\end{sidewaysfigure}

\begin{sidewaysfigure}   
  \begin{center}
    \includegraphics[width=1.1\linewidth] {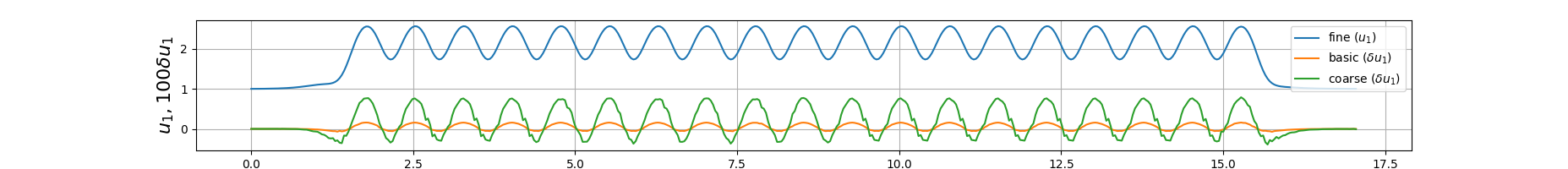} \\
	a \\
    \includegraphics[width=1.1\linewidth] {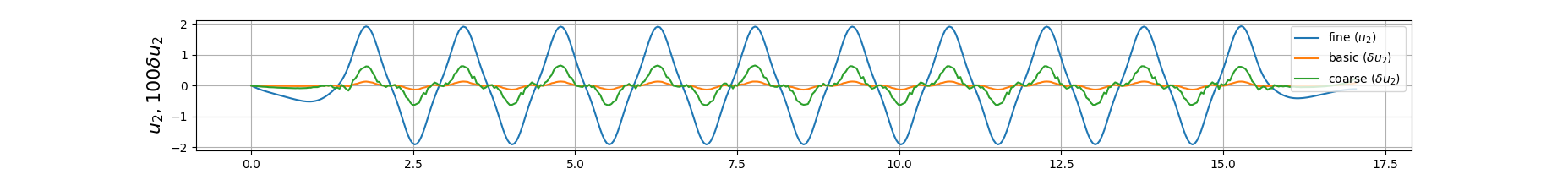} \\
	b \\
    \includegraphics[width=1.1\linewidth] {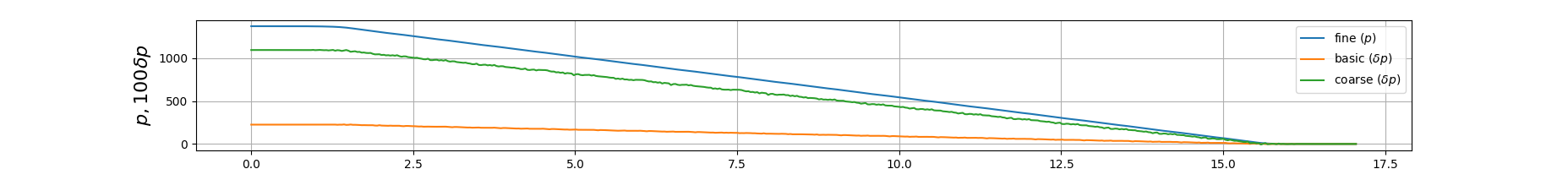} \\
	c \\
  \caption{Velocity components and pressure along the middle line of of the computational domain, computed on fine grid, 
     as well as differences with the solution on coarser grids: 
     a --- velocity component $u_1$, b --- velocity component $u_2$, c --- pressure.   } 
  \label{f-3}
  \end{center}
\end{sidewaysfigure}

As mentioned above, mainly slow flows are of interest for the current study, therefore the basic considered variant is characterized by $\mathrm{Re} =1$. Computed velocity components and pressure are shown on Fig.\ref{f-2}. Convergence of the solution with respect to refinement of the grid is illustrated on Fig.\ref{f-3}. 
We have used three computational grids: basic grid with 18743 nodes and 35958 triangles, coarse grid with 4760  nodes and 8754 triangles and fine grid with 72745  nodes and 142460 triangles.
Fig.\ref{f-3} shows velocity components and pressure on the middle line of the computational domain, computed on the fine grid. The differences ($\delta u_1, \delta u_2, \delta p$) between solutions computed on the fine grid and solutions computed on the basic and on the coarse grid are presented to illustrate the convergence. For convenience, in the visualization these differences are multiplied by 100. From the results, it can be concluded that the basic grid provides good accuracy for the numerical solution.

\subsection{Simulation of reactive transport} 

The unsteady species transport problem (\ref{20}),  (\ref{15}), (\ref{21})-- (\ref{23}) is solved numerically using standard Lagrangian $P_1$ finite elements. Let us define 
\[
 S = \{ s \in H^1(\Omega_f) : \ q = 1 \ \mathrm{on} \  \Gamma_{in} \} ,
\]  
\[
 \hat{S} = \{ s \in H^1(\Omega_f) : \ q = 0 \ \mathrm{on} \  \Gamma_{in} \} .
\] 
The approximate solution $c \in S$ is sought from
\begin{equation}\label{28}
\left (\frac{\partial c}{\partial t}, s \right) + d(c,s) =  
 \left (\frac{\partial m}{\partial t} ,s \right )_s
 \quad \forall s \in \hat{S} ,   
\end{equation} 
where the following notations are used
\[
  d(c,s) := - \int_{\Omega_f} c \bm u \cdot \nabla  s \, d \bm x 
  + \frac{1}{\mathrm{Pe}} \int_{\Omega_f} \nabla c \cdot \nabla  s \, d \bm x 
  + \int_{\Gamma_{out}} (\bm u \cdot \bm n) c s \, d \bm x ,
\] 
\[
  (\varphi,s)_s := - \int_{\Gamma_{s}} \varphi   s \, d \bm x .
\] 
For determining $m \in G = L_2(\Gamma_{s})$ (see (\ref{11})) we use
\begin{equation}\label{29}
 \left (\frac{\partial m}{\partial t} , g \right )_s -
 (f(c,m), g)_s = 0,
 \quad g \in G .  
\end{equation} 
The discretization in time is based on symmetric discretization (Crank–Nicolson method), which is second order accurate (see, e.g., \cite{Samarskii,Ascher2008}. 
Let $\tau$  be a step-size of a uniform grid in time such that 
$c^n = c(t^n), t^n = n\tau, \ n = 0, 1, ...$. 

Eq.(\ref{28}) is approximated in time as follows
\[
\left (\frac{c^{n+1} - c^n}{\tau }, s \right) + d\left (\frac{c^{n+1} + c^n}{2} ,s\right ) =  
 \left (\frac{m^{n+1} - m^n}{\tau } ,s \right )_s .
\] 
Similarly, for (\ref{29}) we get
\[
\left (\frac{m^{n+1} - m^n}{\tau }  , g \right )_s -
 \left (f\left (\frac{c^{n+1} + c^n}{2} ,\frac{m^{n+1} + m^n}{2} \right ), g\right )_s = 0 ,
 \quad n = 0, 1, ... , 
\] 

In the considered here case, zero initial conditions are posed
\[
 c^0 = 0, \quad \bm x \in \Omega_f, 
\] 
\[
 m^0 = 0, \quad \bm x \in \Gamma_s . 
\]  

Henry adsorption isotherm (\ref{24}) is used in the simulations reported below. The basic set of parameters reads as follows:
\[
 \mathrm{Pe} = 10,
 \quad \mathrm{Da}_a = 0.005,
 \quad \mathrm{Da}_d = 0.05,  
\] 
The unsteady problem is solved on dimensionless time interval $(0,T)$, $T = 40$, using time step $\tau =0.1$. 
Breakthrough curve (i.e., average outlet concentration) is used to characterize the reactive transport:
\begin{equation}\label{30}
 c_{out}(t) = \frac{\int_{\Gamma_{out}} c(\bm x, t) d \bm x}{\int_{\Gamma_{out}} d \bm x} . 
\end{equation} 

\begin{sidewaysfigure}   
  \begin{center}
    \includegraphics[width=0.5\linewidth] {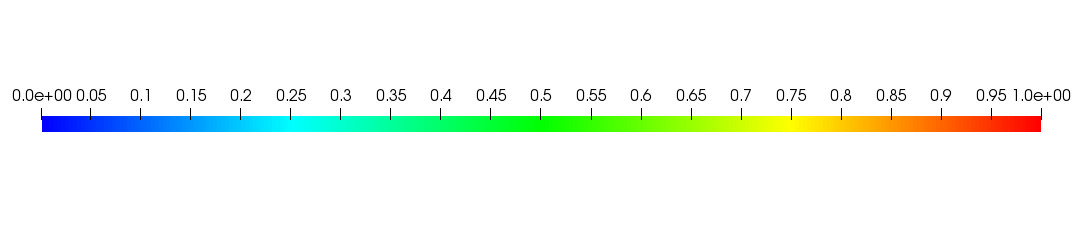} \vspace{-8mm} \\
    \includegraphics[width=1.1\linewidth] {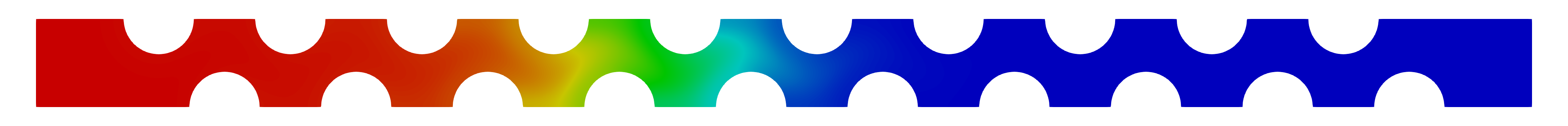} \\
	 a \\
    \includegraphics[width=1.1\linewidth] {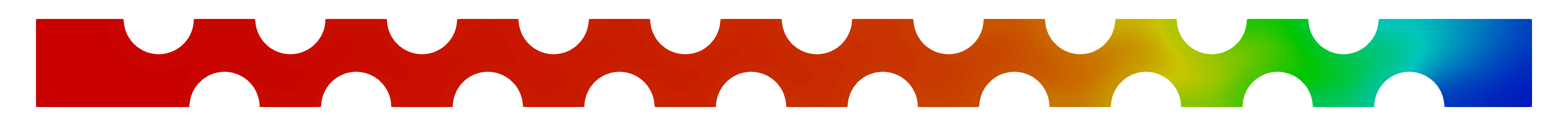} \\
	 b \\
    \includegraphics[width=1.1\linewidth] {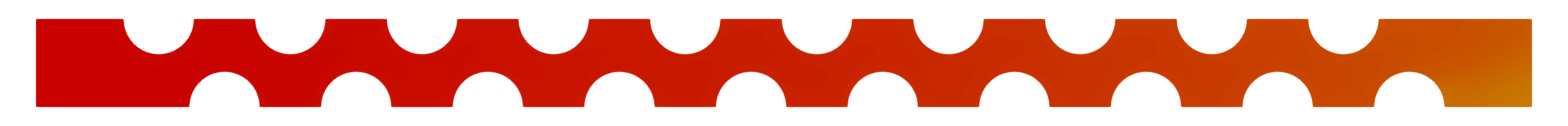} \\
	 c 
  \caption{Concentration $c$ at different time moments, a --- $t = 5$,	b --- $t = 10$, c --- $t = 15$}
  \label{f-4}
  \end{center}
\end{sidewaysfigure}

\begin{sidewaysfigure}   
  \begin{center}
    \includegraphics[width=1.1\linewidth] {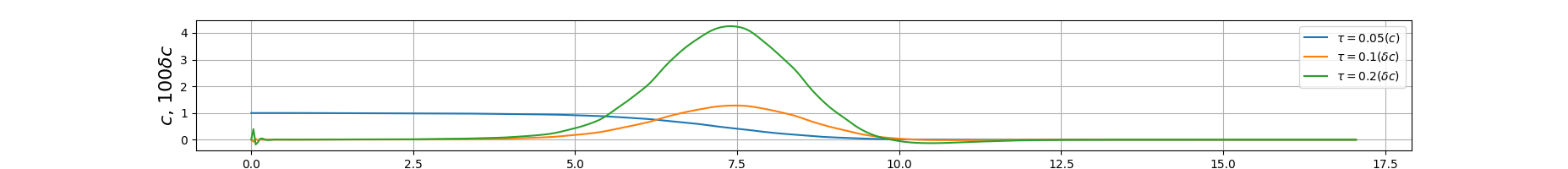} \\
	a \\
    \includegraphics[width=1.1\linewidth] {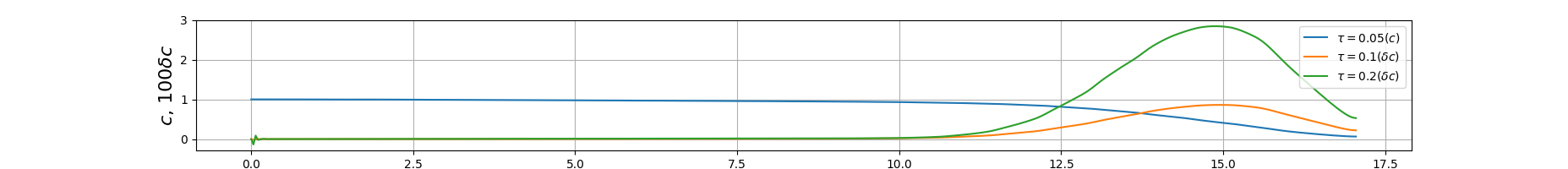} \\
	b \\
    \includegraphics[width=1.1\linewidth] {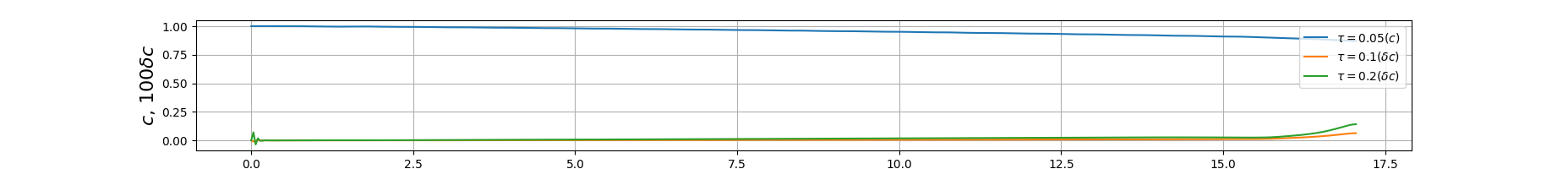} \\
	c \\
  \caption{Concentration $c$ computed with $\tau = 0.05$ and its difference $\delta c$ 
           from the concentrations computed with $\tau = 0.1$   and $\tau = 0.2$ at different time moments: 
           a --- $t = 5$, b --- $t = 10$, c --- $t = 15$}
  \label{f-5}
  \end{center}
\end{sidewaysfigure}

Species concentration at different time moments are shown on Fig.\ref{f-4}. 
The accuracy of the approximate solution with respect to the time step is illustrated on Fig.\ref{f-5}. 
The solution along the middle line of the computational domain is shown there for the case $\tau = 0.05$, 
as well as the magnified by factor 100 differences between this solution and solutions obtained 
with $\tau = 0.1$ and $\tau = 0.2$. 
The results show that a good accuracy can be achieved using $\tau = 0.1$, 
and the further simulations are performed using this time step.

\begin{figure}[htp]
  \begin{center}
 	\includegraphics[width=0.49\linewidth] {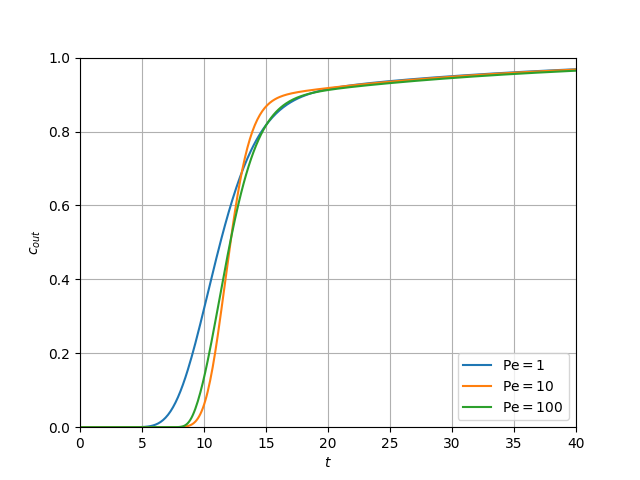} 
 	\includegraphics[width=0.49\linewidth] {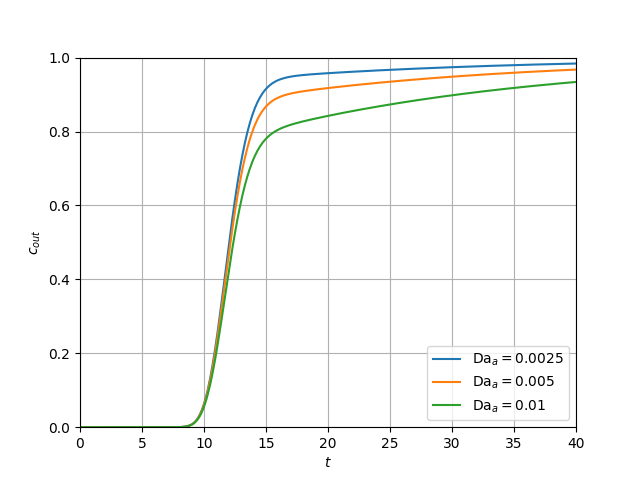} \\
	\makebox[0.49\linewidth][c] {a} \makebox[0.49\linewidth][c] {b}  
 	\includegraphics[width=0.49\linewidth] {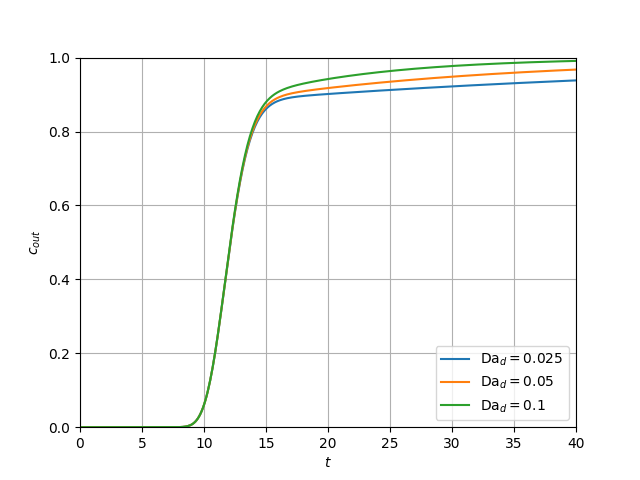} 
 	\includegraphics[width=0.49\linewidth] {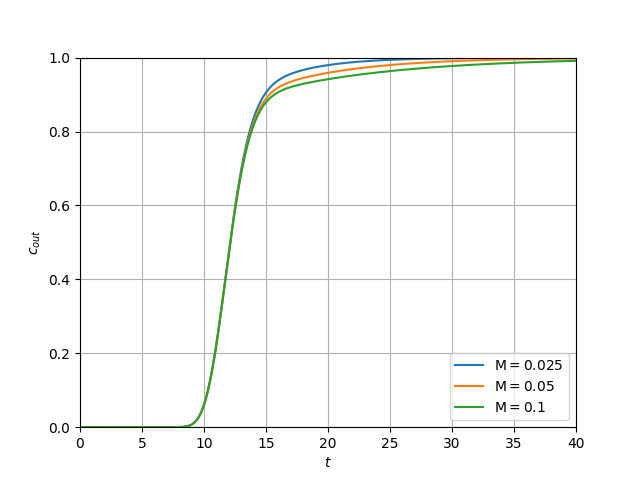} \\
	\makebox[0.49\linewidth][c] {c} \makebox[0.49\linewidth][c] {d}  
	\caption{Outlet concentration: a --- dependence from  $\mathrm{Pe}$,
            b --- dependence from $\mathrm{Da}_a$, c --- dependence from $\mathrm{Da}_d$, d --- dependence from $\mathrm{M}$}
	\label{f-6}
  \end{center}
\end{figure}

Sensitivity study is carried out in order to see how the change of different parameters leads to change of the breakthrough curves. Such sensitivity studies are often first stage for optimization or parameter identification procedures. The dependence of the average outlet concentration from Peclet number is shown on Fig.\ref{f-6}a. Smaller $\mathrm{Pe}$ means larger diffusion, what is well illustrated by the spreading of the outlet concentration.
Note that the amount of the deposited concentration is given by $1- c_{out}(t)$.

The rate of adsorption is characterized by $\mathrm{Da}_a$. The influence of this parameter on the outlet concentration is illustrated on Fig.\ref{f-6}b. Increasing $\mathrm{Da}_a$, as expected, leads to more intensive adsorption and larger amount of the deposited substance. The influence of the parameter $\mathrm{Da}_d$ on the outflow concentration is illustrated on Fig.\ref{f-6}c. In the current paper we restrict our considerations to the case of small Pe and Damkoler numbers, the study for other regimes will be presented in forthcoming papers. 

Henry isotherm usually describes well the initial stages of the adsorption. In cases when only a limited mass can be adsorbed at a surface, Langmuir isotherm should be used to reflect the decay of the adsorption rate close to the saturation. In this case an additional parameter appears, namely $\mathrm{M}$ (see (\ref{25})).
The influence of $\mathrm{M}$ on the average output concentration is shown on Fig.\ref{f-6}d.

\clearpage
\section{Numerical solution of the inverse problem} 

Consider an inverse problem for determining unknown adsorption rate (so called parameter identification problem), using information about the dynamics of the average outlet concentration $\widetilde{c}_{out}(t)$ (the latter is usually easy to measure in experiments). The starting point is monitoring of the difference (residual) between the measured $\widetilde{c}_{out}(t)$ and the computed  $c_{out}(t)$ average outflow concentration for different values of the parameters $\mathrm{Da}_a$ and $\mathrm{Da}_d$ in Henry isotherm Eq.(\ref{24}). 

\subsection{Direct computation of the residual functional} 

Let us consider the problem for identifying the adsorption and desorption rates in Langmuir isotherm, $\mathrm{Da}_a, \mathrm{Da}_d$, respectively, using given measurement results $\widetilde{c}_{out}(t)$. The functional of the residual is given by
\begin{equation}\label{31}
 J(\mathrm{Da}_a, \mathrm{Da}_d) = \int_{0}^{T} (c_{out}(t) - \widetilde{c}(t))^2 d t,
\end{equation} 
where $c_{out}(t)$ is given by (\ref{30}) after solving the dimensionless form of (\ref{3}) together with the respective boundary conditions, and after solving (\ref{20})--(\ref{23}). We use the so called synthetic approach to emulate measurements data. At a preliminary stage, the direct problem is solved for  selected fixed values of ($\widetilde{\mathrm{Da}}_a, \ \widetilde{\mathrm{Da}}_d$) (called later exact solution of the identification procedure, or exact parameters, here $(0.005,0.05)$),and the computed $c_{out}(t)$ is used as measurement data in the case of no noise, i.e., in this case $\widetilde{c}(t)=c_{out}(t)$. 

\begin{figure}[htp]
  \begin{center}
 	\includegraphics[width=0.75\linewidth] {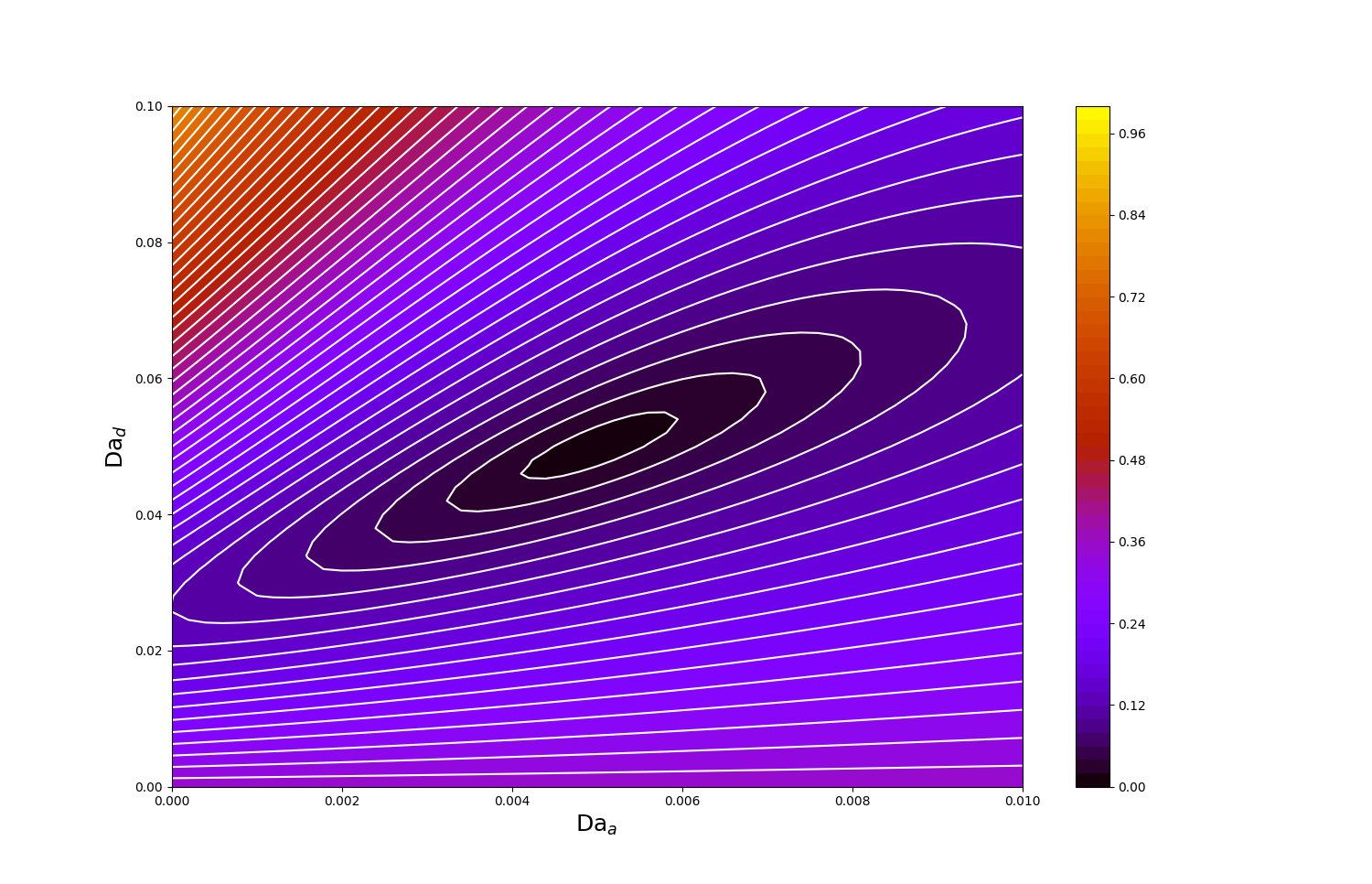} \\
	a \\
  	\includegraphics[width=0.75\linewidth] {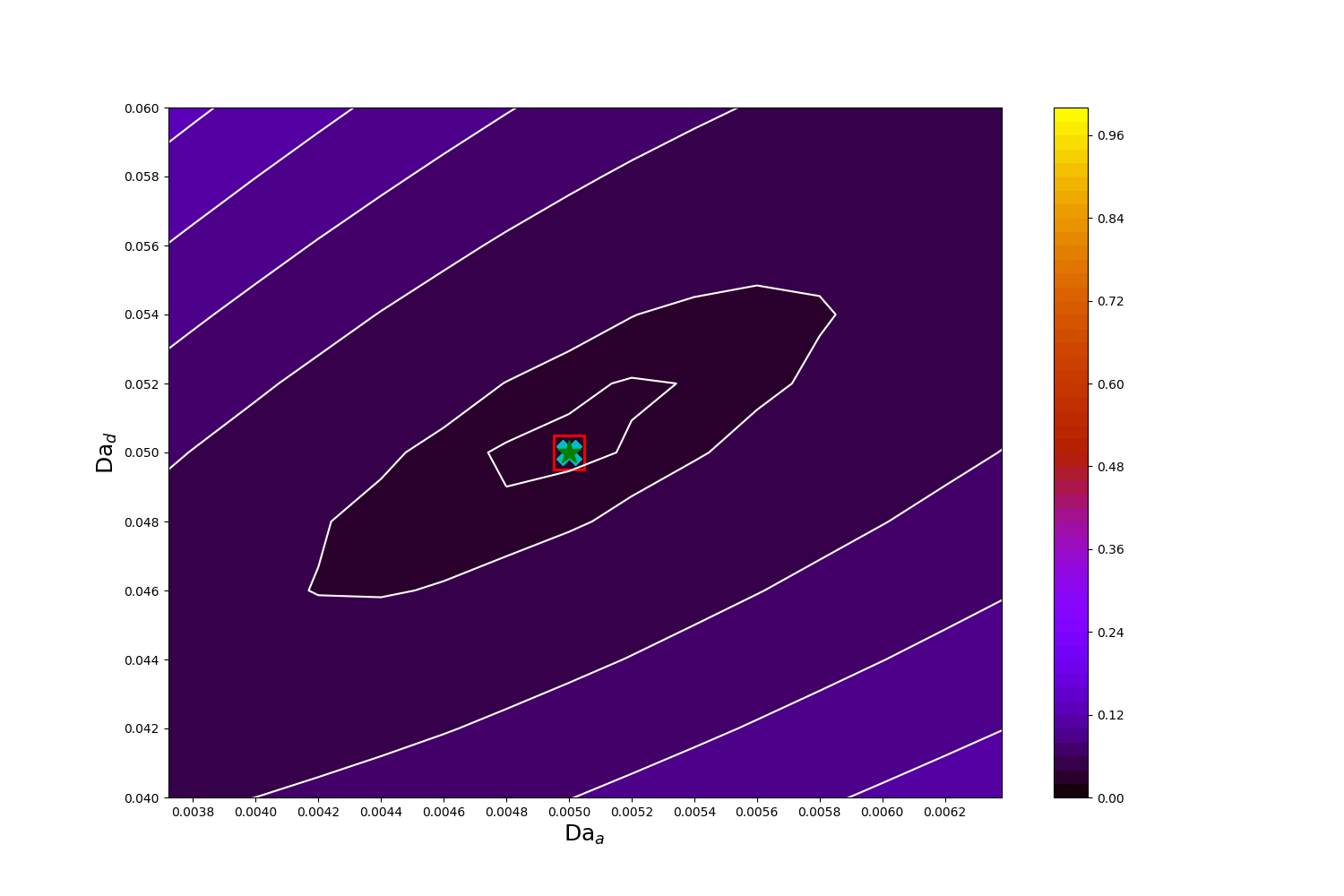}  \\
	b \\
 	\includegraphics[width=0.75\linewidth] {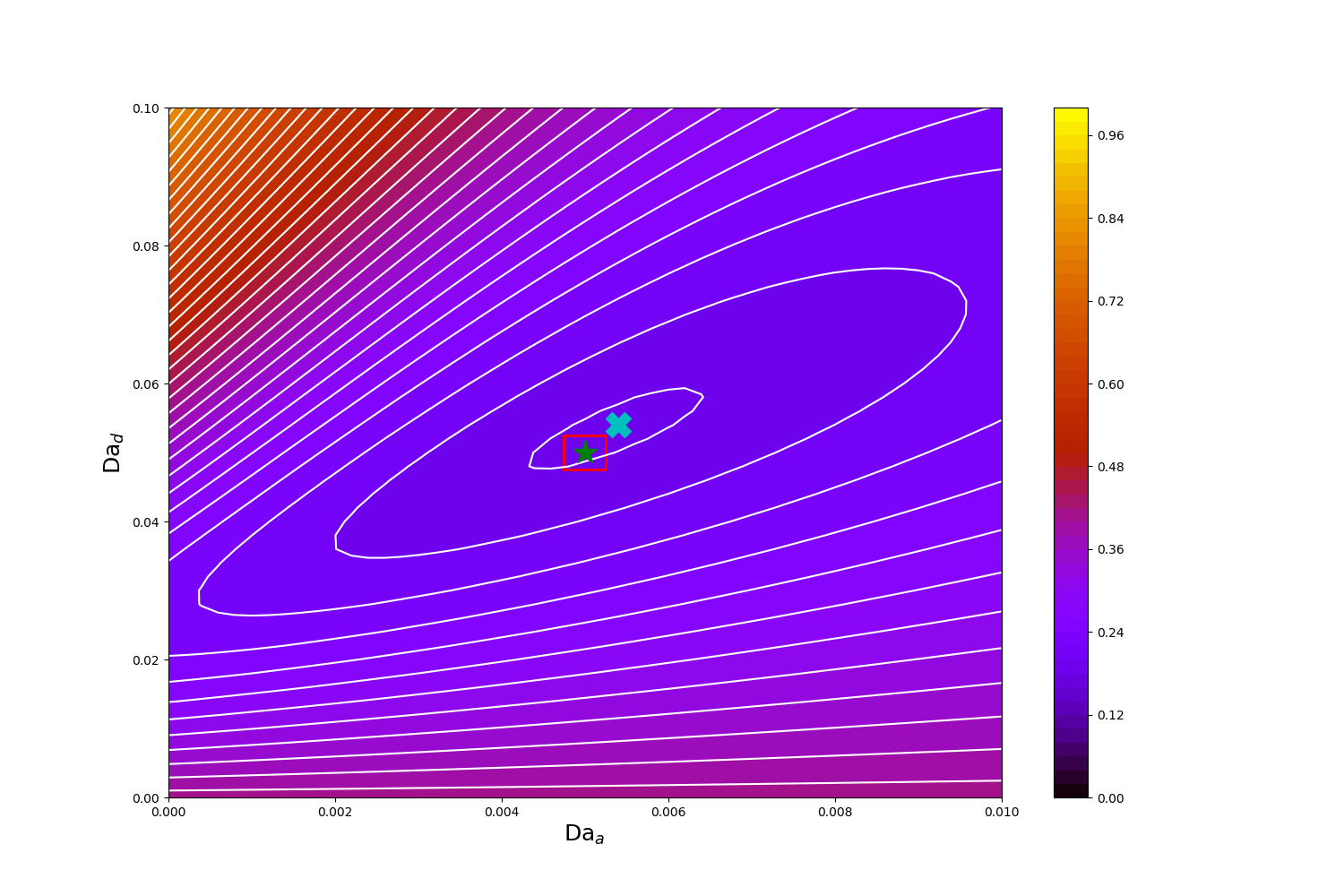} \\
	c \\
	\caption{Square root of the residual functional ($J^{1/2}(\mathrm{Da}_a, \mathrm{Da}_d)$): a --- exact data ($\delta = 0$),
                  b --- $\delta = 0.01$,  c --- $\delta = 0.05$}
	\label{f-7}
  \end{center}
\end{figure}

First of all, having in mind that this is a relatively simple parameter identification problem (only two parameters need to be identified, the direct problem can be fast and efficiently solved numerically), a simple strategy for parameter identification will be considered here. Consider starting feasible set $G$ defined as
\[
 0 \leq \mathrm{Da}_a \leq 0.01, 
 \quad 0 \leq  \mathrm{Da}_d \leq 0.1 .
\]  
Computations are performed on uniform grid for the parameters $\mathrm{Da}_a, \ \mathrm{Da}_d$ with $51\times 51$ nodes, and the direct problem is solved for each pair of parameters from this grid. Square root from the residual functional 
is shown on Fig.\ref{f-7}a. Here and below the isolines are visualized using the library {\it matplotlib}. The isolines of the functional are drawn with step $0.02$. In this case the minimum of the functional is zero, $\min J(\mathrm{Da}_a, \mathrm{Da}_d) = 0$, and it is achieved at point $(0.005,0.05)$). Recall that this point was used to emulate by computations the measurement data.  In this particular case this point belongs to the parameter grid on which we seek the solution, therefore $\min J(\mathrm{Da}_a, \mathrm{Da}_d) = 0$ in this case.

\subsection{The influence of the noise in the measured data} 

In the reality the measurement data are often noisy, therefore it is important to understand how the noise influences the parameter identification. Suppose that the measurement are carried out at each time step. Instead of (\ref{30}), let us consider noisy measurement:
\[
  \widetilde{c}(t^n) = c_{out}(t^n; \widetilde{\mathrm{Da}}_a, \widetilde{\mathrm{Da}}_d) +
  \delta \sigma(t^n),
  \quad n = 1,2, ..., N, 
  \quad N \tau = T .   
\] 
Here the parameter $\delta$ quantifies the amplitude of the noise, while $\sigma(t^n)$ is a random variable, uniformly distributed on the interval $[-1,1]$.

In the case of noisy measurement data, instead of trying to identify one point in the parameters space (what is the most often considered case in parameter identification problems), one should identify a set in the parameters space, for which the value of the residual functional is below prescribed tolerance. Having in mind the dispersion of the uniformly distributed random variable $\sigma$, we define the set of the admissible parameters by the inequality  
\begin{equation}\label{33}
J(\mathrm{Da}_a, \mathrm{Da}_d) \leq  \gamma \delta^2 \frac{T}{3} , 
\end{equation} 
where $\gamma > 1$ is a numerical parameter.
Here and below in order to avoid misunderstanding, we will call {\it admissible set, or resultimg admissible set, or final admissible set} the set which is identified as a result of the identification, thus making difference with {\it feasible set or starting feasible set}, which is the parameter set $G$ which is input for our identification procedure. 

\begin{figure}[htp]
  \begin{center}
 	\includegraphics[width=0.49\linewidth] {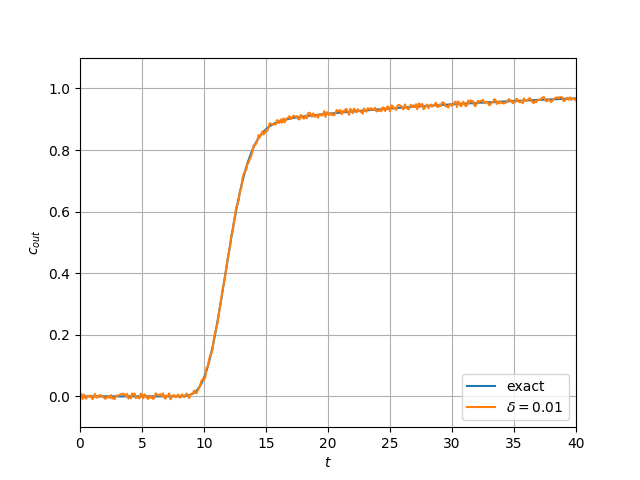} 
 	\includegraphics[width=0.49\linewidth] {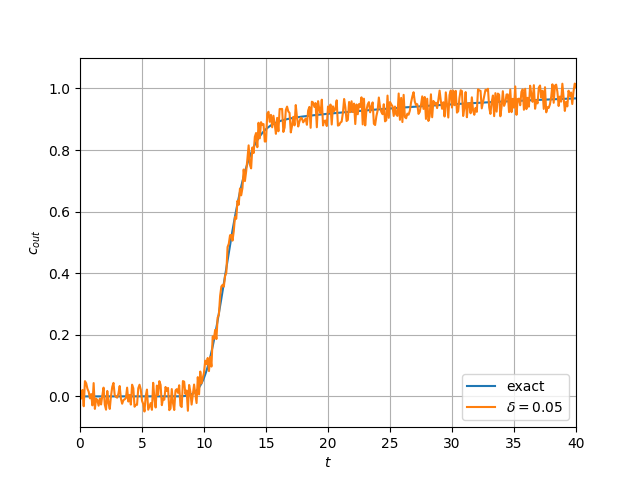} \\
	\makebox[0.49\linewidth][c] {a} \makebox[0.49\linewidth][c] {b}  
	\caption{Average outflow concentration: a --- $\delta = 0.01$, b --- $\delta = 0.05$}
	\label{f-8}
  \end{center}
\end{figure}

In the above a priori estimate the definition of the admissible set explicitly depends on the time interval during which the measurements are done. Obviously, it may be reasonable to ignore the measurement data before the breakthrough, as well as to ignore the measurement data after an equilibrium is achieved. The question about the minimum required information for an identification procedure, e.g., with which time interval it is reasonable to work, is not discussed here. Relatively large time interval for the measurements is considered here, see Fig.\ref{f-8}. 

One realization of the noisy data for the the average outlet concentration is shown on Fig.\ref{f-8}a for the case $\delta = 0.01$. Square root of the residual functional for this case is visualized by color on Fig.\ref{f-7}b, some isolines are also plotted. In this case, due to the noise, the minimum of the functional is not equal to zero, it is equal to $0.0361$. The point at which the functional takes its minimum is marked with $\times$ on the figure. The exact parameters are marked with $\star$ on the figure.
The red square on the picture, centered around the star (which almost coincides with the $\times$ in this case), is a drawing of $(0.99\times1.01)\mathrm{Da}_a,  (0.99\times1.01)\mathrm{Da}_d$ isolines of the residual. That is, within this square the error in identifying the exact parameters is less than $1\%$. Furthermore, the first plotted isoline for the residual $J$, is computed from (\ref{33}) for $\gamma=1.02625$.

Similar results for significantly higher amplitude of the noise, ($\delta = 0.05$, Fig.\ref{f-8}b) are shown on Fig.\ref{f-7}c. In this case, the minimum of the functional is equal to $0.1840$. The red square in this case indicates $5\%$ error, to be in the same range as the noise. In this case the first plotted isoline for the residual $J$ is also calculated with $\gamma=1.02625$. We see that the minimum of the functional in this case is achieved in a point which is relatively far from the exact value, indicating that one measurement can not be enough to identify the parameters accurately in the case of very noisy data. This is manifestation of the known fact that more Monte Carlo simulations (simulations with different noise realizations, see the next paragraph) have to be done for larger amplitude of the noise. 

\begin{figure}[htp]
  \begin{center}
    \includegraphics[width=1\linewidth]{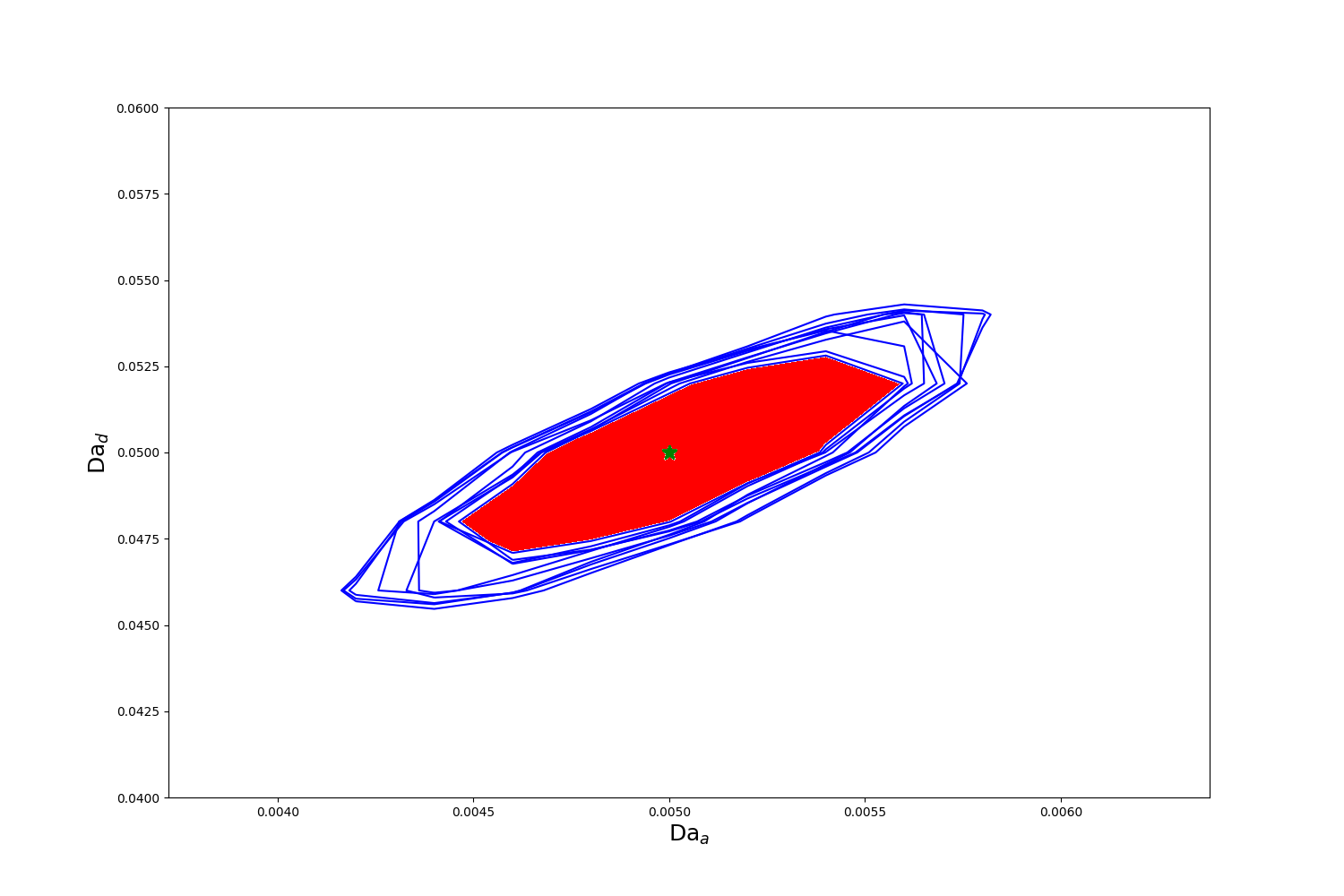}
	\caption{Admissible set of adsorption parameters for ten realizations for noise amplitude $\delta = 0.01$}
	\label{adm001}
  \end{center}
\end{figure}

\begin{figure}[htp]
  \begin{center}
    \includegraphics[width=1\linewidth]{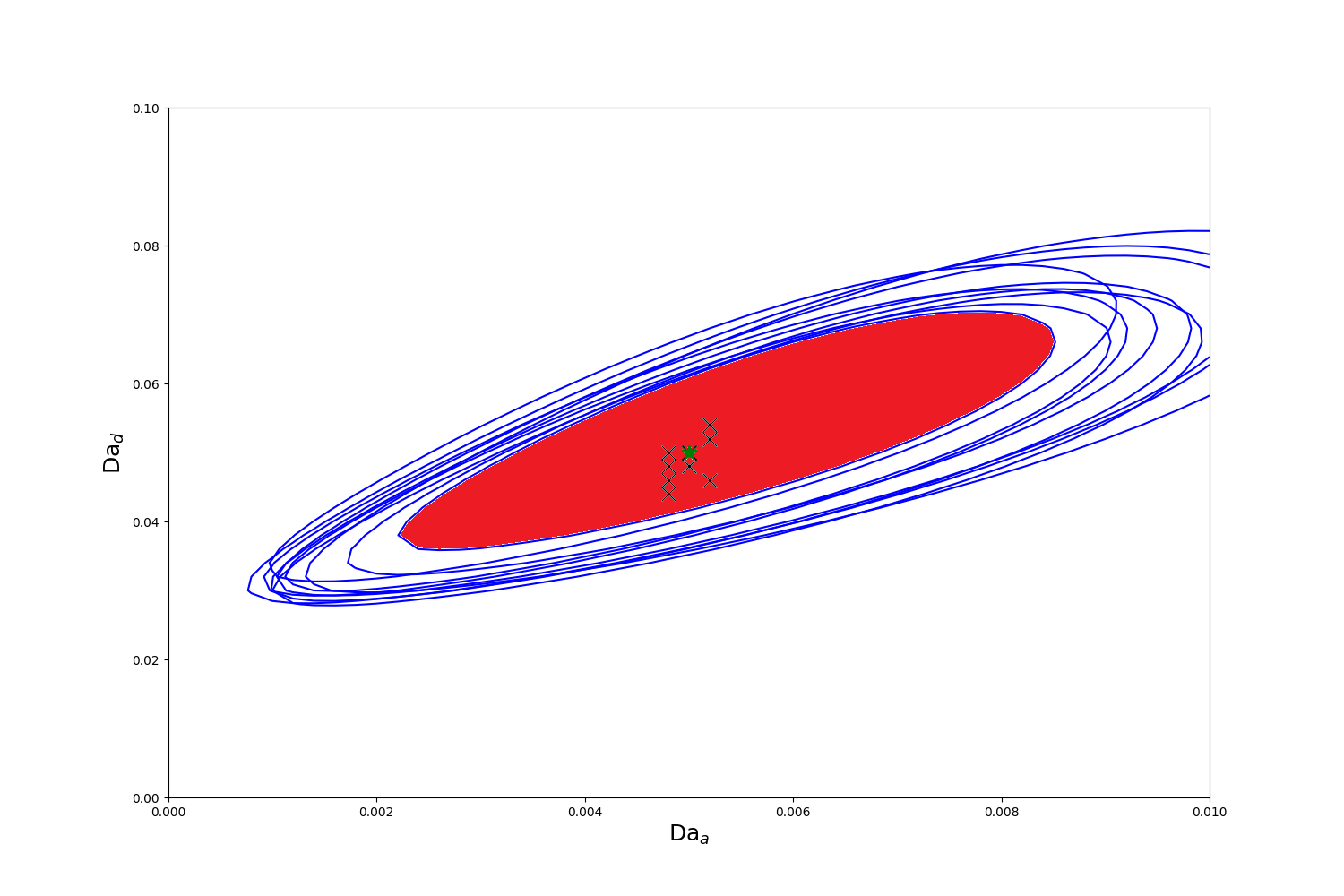}
	\caption{Admissible set of adsorption parameters for ten realizations for noise amplitude $\delta = 0.05$}
	\label{adm005}
  \end{center}
\end{figure}

To investigate the role of the stochasticity in the noise, simulations with ten different realizations of the noise are done for each of the noise levels $\delta = 0.01$ and $\delta = 0.05$. The resulted admissible sets are shown on Fig.\ref{adm001} and Fig.\ref{adm005}. For a better visualization, here one isoline is plotted for each realization, and relatively large $\gamma=1.21$ is selected for this isoline. The points in which the minima of the functional are achieved are close to the exact point for $\delta=0.01$ and therefore they are not plotted. For $\delta=0.05$ they are denoted with $\times$ on the figure. The overlapping area is shrinking, indicating that the accuracy of the parameter identification can be improved by performing a number of simulations for different realizations of the noise. It should be noted that in \cite{maday2015pbdw} the authors have investigated the role of the noisy data in connection with another parameter identification problem. They have proven that under certain conditions on the noise (which are satisfied also in our case), within a Monte Carlo procedure (working with sufficiently large set of measurement data), the exact parameters can be identified. However, very large number of measurements may be needed in the case of large variance. Such an approach is important in the case when the parameters have to be identified with high accuracy. However, in many industrial problems the goal is to identify the parameters with an accuracy which is of the order of the noise. For such problems, it makes sense to identify admissible sets of parameters, instead of trying to identify the exact parameters. 

It should be noted that it is not granted that the computed final admissible sets for two different realizations of the noise will overlap. If this will be the case with some of the computed final admissible sets, a subset of overlapping admissible sets has to be chosen in order to select the shrinked area. We leave the discussion on Monte Carlo and other approaches for more accurate parameter identification in the case of noisy data for a future article.

\subsection{Statistical parameter identification} 

The problem for identifying adsorption parameters in this case is formulated as follows.
Find resulting (final) admissible set $\mathrm{Da}_a, \mathrm{Da}_d$ satisfying (\ref{33}) from a starting set of admissible parameters $G$. 
This problem can be solved by different local, and in general case, global optimization approaches \cite{horst2013handbook,nocedal2006numerical}. 

Recall, that the standard approach for solving parameter identification problems for PDEs \cite{samarskii2007numerical,kern2016numerical} does not target at determining an admissible set of parameters, instead, it targets at determining a particular point in the parameter space. Usually gradient methods are used as a part of such approaches. The evaluation of the admissible set of parameters for such deterministic approaches can be based on iterative minimization of the residual functional starting from different initial values, the so called multistart method.  

As alternative we consider here stochastic methods and their variants based on deterministic sampling of points. A stochastic approach for global optimization in its simplest
form consists only of a random search and it is called Pure Random Search 
\cite{Zhigljavsky2008}. In this case the residual functional 
is evaluated at randomly chosen points belonging to the feasible (initial) parameter set. Such an stochastic approach is more efficient than the computation of the residual function on uniform grid in the feasible set. 

Consider the case when the starting feasible set $G$ is defined in a simple way:
\[
 0 \leq \mathrm{Da}_a \leq \overline{\mathrm{Da}}_a,
 \quad  0 \leq \mathrm{Da}_d \leq \overline{\mathrm{Da}}_d .
\] 
To generate of uniformly distributed points $D$ from $G$, Sobol sequences \cite{sobol1976uniformly,sobol1979systematic} are used. The software implementation is based on the sensitivity analysis library SALib \cite{Herman2017}. Such an approach is widely used in applied problems for multicriteria parameter identification (see., e.g. \cite{sobol1981choosing}).

\begin{figure}[htp]
  \begin{center}
  	\includegraphics[width=0.75\linewidth] {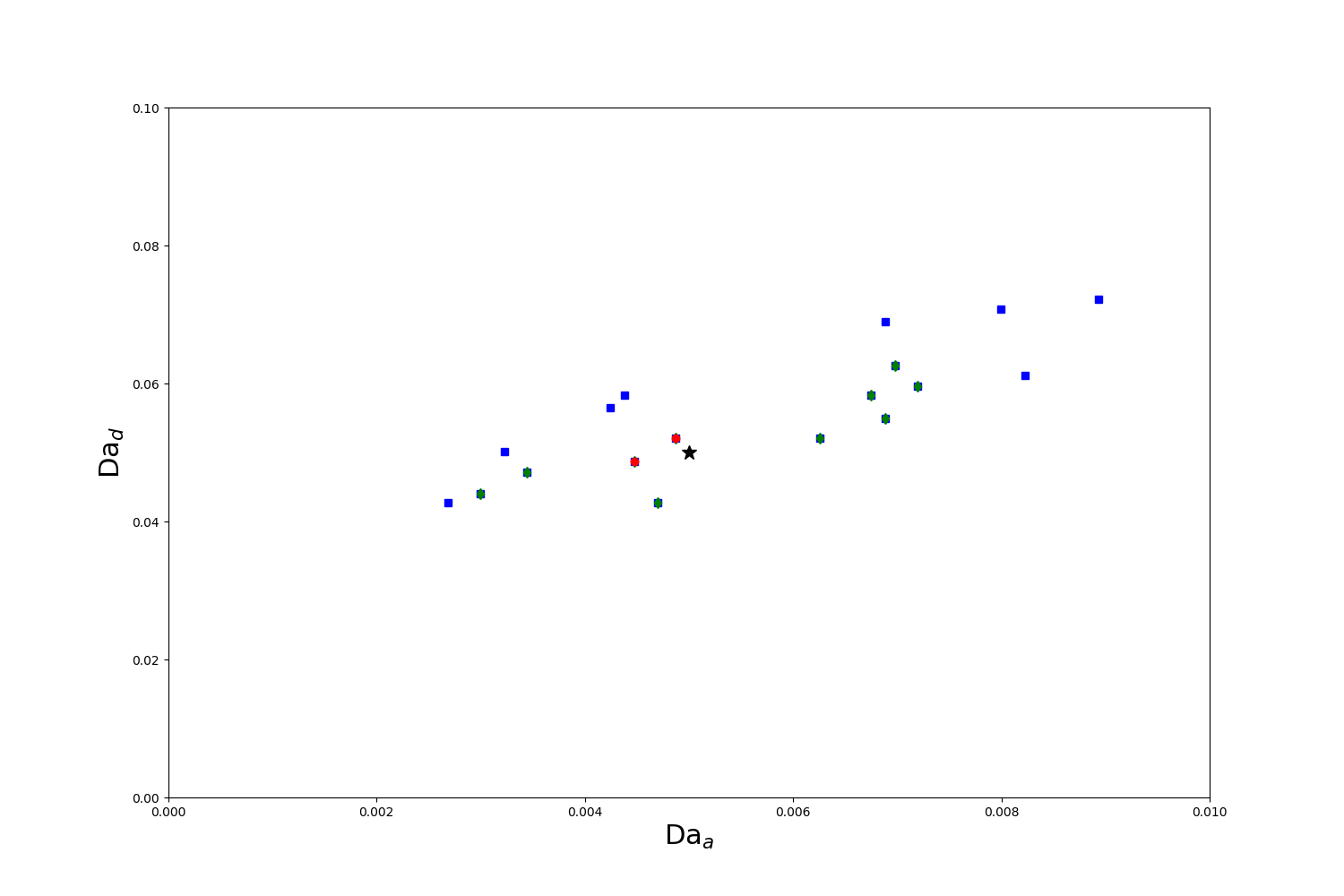} \\
	a \\
  	 	\includegraphics[width=0.75\linewidth] {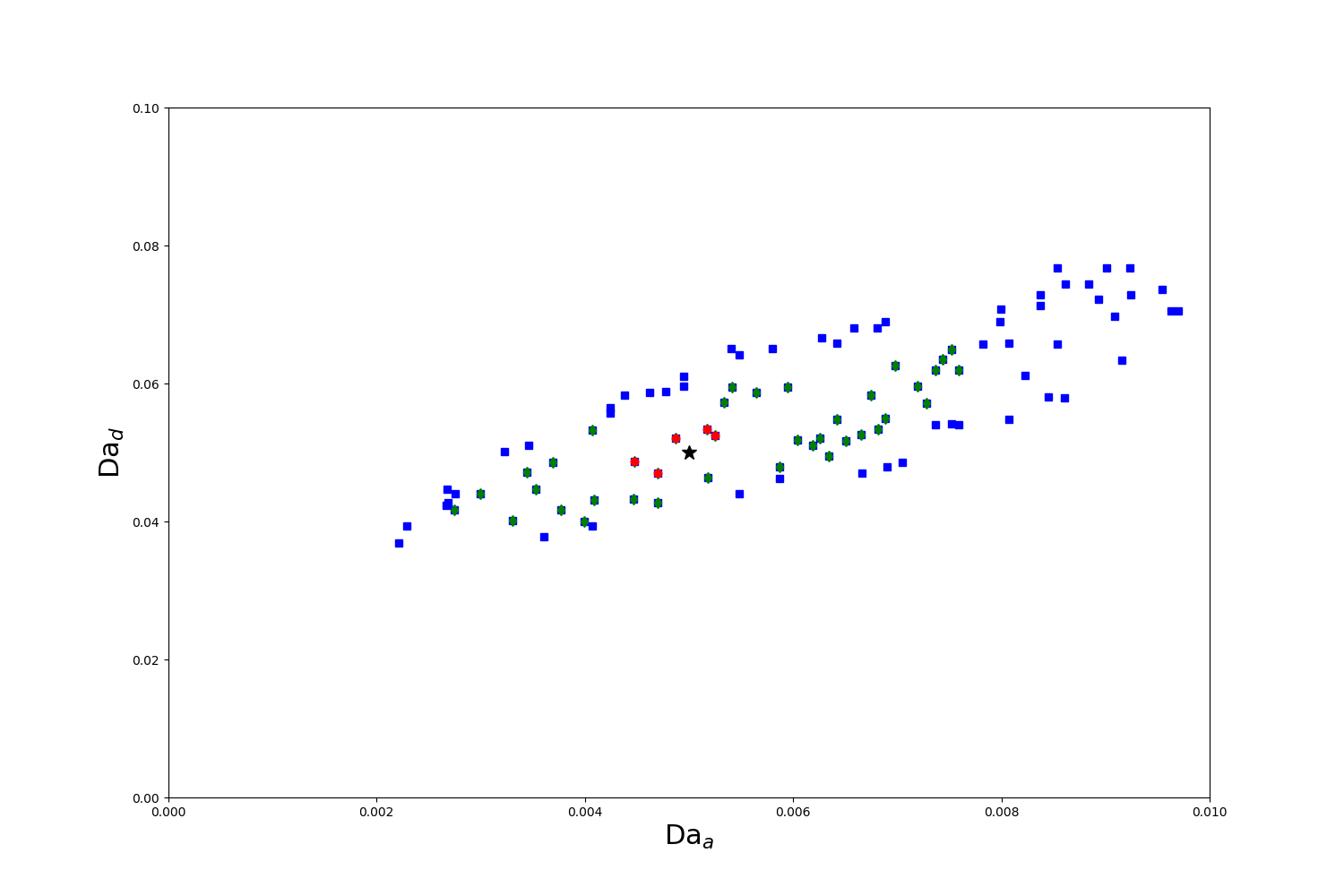} \\
	b \\
	\caption{Deterministic sampling of points, plot of points for which the inequality (\ref{33}) is satisfied, $\gamma=1.02625$ ($\delta = 0.01$ --- red, $0.03$ --- green,  $0.05$ --- blue): a --- $N = 150$,
 	b --- 
 $N = 600$}
	\label{f-10}
  \end{center}
\end{figure} 

The result of the statistical identification for 
$N = 150$ sampling points is shown on Fig.\ref{f-10}a. Similar results for 
$N = 600$ sampling points is shown on Fig.\ref{f-10}b 
In both cases we are plotting the points from the random set for which the value of the residual satisfies (\ref{33}) with $\gamma=1.02625$. 
It can be seen that for small number of the sampling points it is difficult to reconstruct reasonable final admissible set because only few points satisfy (\ref{33}). More points are available for larger amplitude of the noise because of the flattening of the functional, but at the same time the accuracy is worse in this case. 

\subsection{Multistage parameter identification} 

\begin{figure}[htp]
  \begin{center}
    \includegraphics[scale = 0.3] {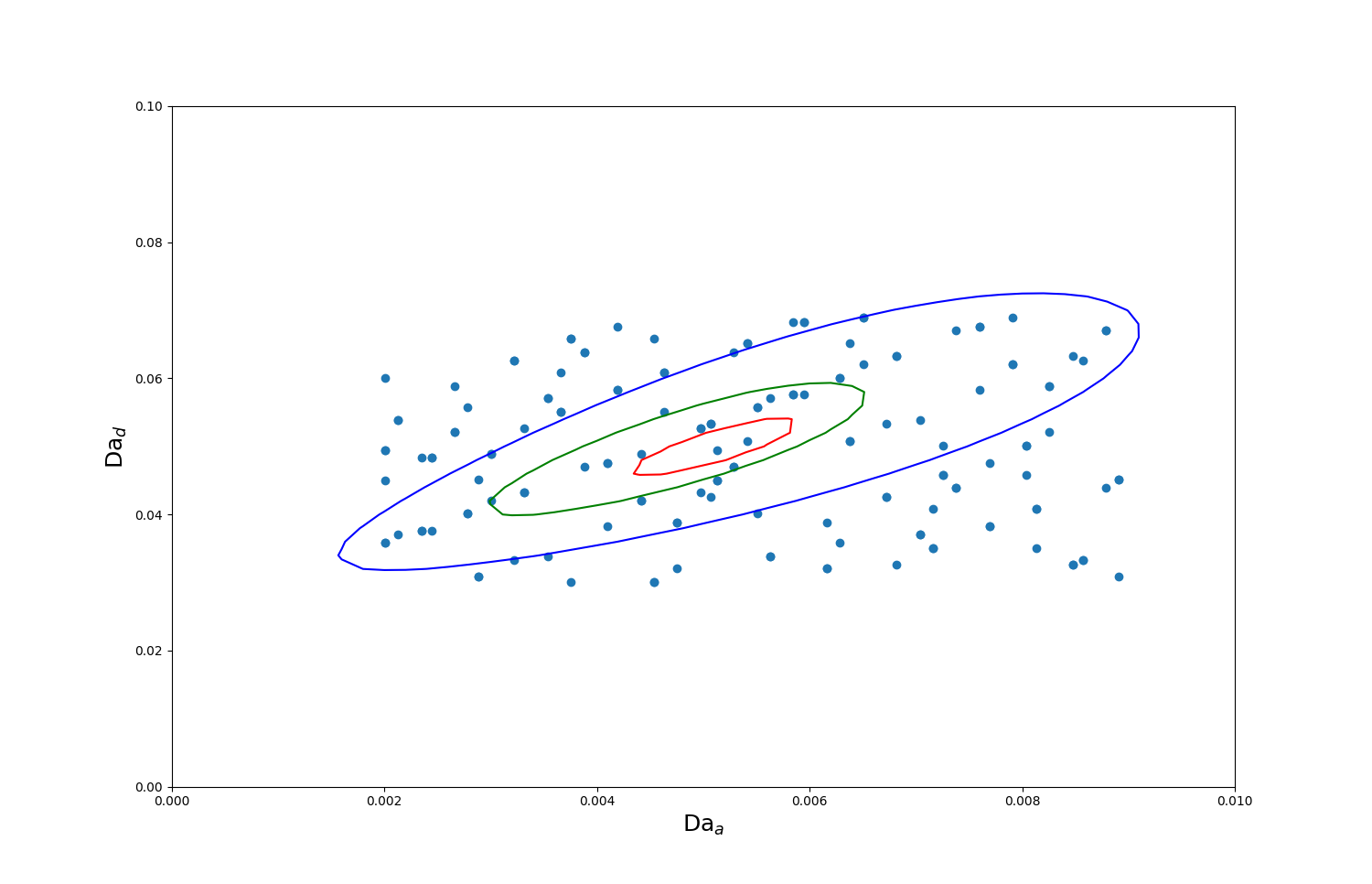}
	\caption{Shrinked domain of feasible parameters for $N = 150$: $\delta = 0.01$ --- red, $0.03$ --- green,  $0.05$ --- blue}
	\label{f-11}
  \end{center}
\end{figure}

In the considered here case, when only two parameters need to be identified and we deal with convex functional, one can consider multistage parameter identification (could also be called predictor-refinement parameter identification).  Instead of increasing the number of the samples, one can consider a multistage method, consecutively shrinking the starting feasible set $G$. For example, the results (see Fig.\ref{f-7}a) obtained with $N = 150$ samples for the original initial feasible set $G$
\[
0. \leq \mathrm{Da}_a \leq 0.01,
\quad  0. \leq \mathrm{Da}_d \leq 0.1 .
\] 
can be evaluated and a new, smaller starting feasible set can be defined for the second stage of the identification procedure, e.g.;
\[
 0.002 \leq \mathrm{Da}_a \leq 0.009,
 \quad  0.03 \leq \mathrm{Da}_d \leq 0.07 .
\] 
The simulation with $N = 150$ sample in the reduced domain at the second stage of the identification procedure, give good accuracy in determining the final admissible set, see Fig.\ref{f-11}. 

\begin{figure}[htp]
  \begin{center}
    \includegraphics[scale = 0.3] {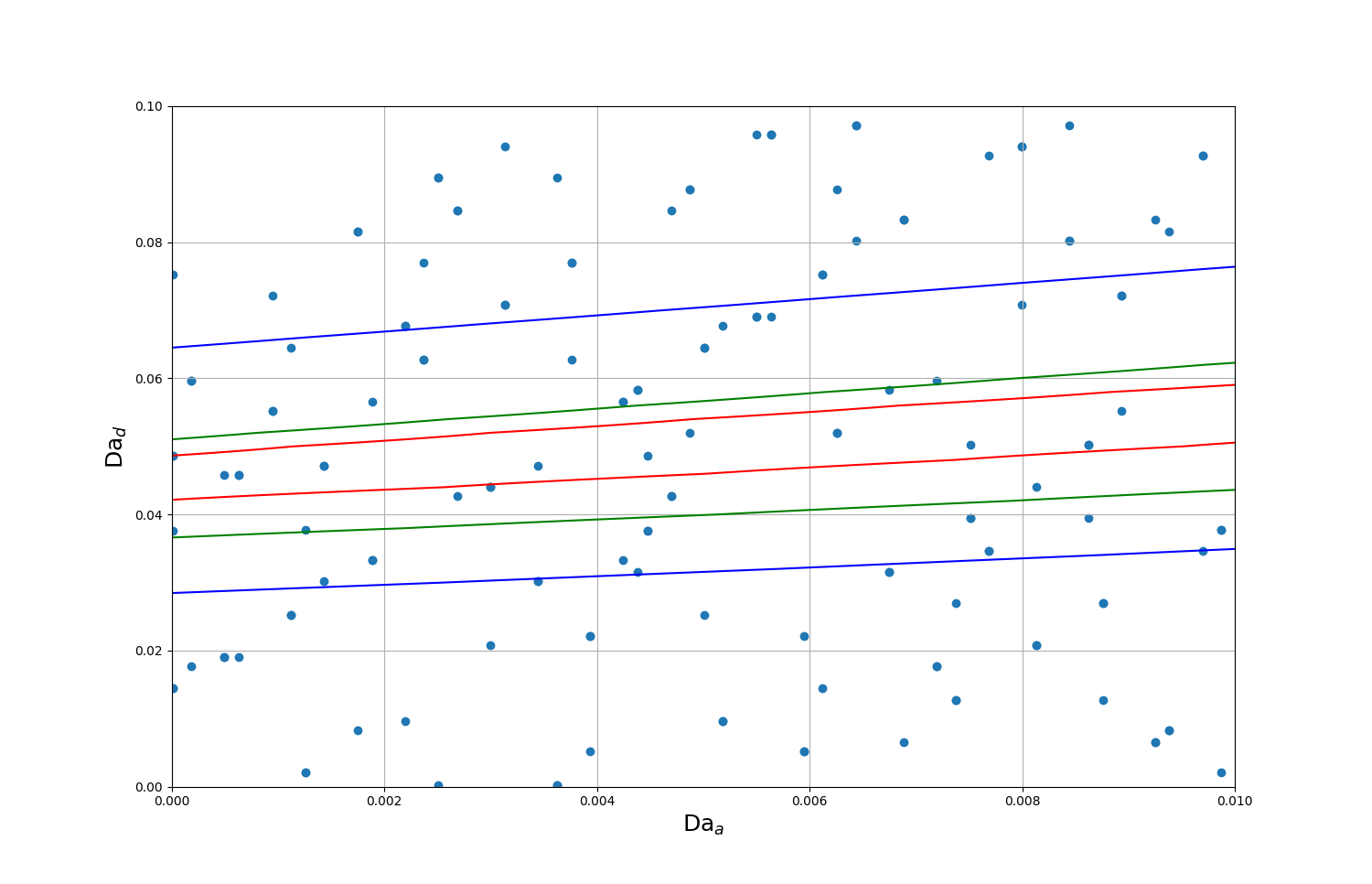}
	\caption{Usage of input data till $T = 15$ with $N = 150$: $\delta = 0.01$ --- red, $0.03$ --- green,  $0.05$ --- blue}
	\label{f-12}
  \end{center}
\end{figure}

Additional peculiarity of the considered here problem is the fact that different adsorption parameters have different influence on the solution at different time intervals. For example, results from Fig.\ref{f-6} indicate that the dynamic of the outflow concentration is different on the subintervals left and right from $T = 15$. This makes reasonable the consideration of the residual functional (\ref{31}) at a smaller time interval, namely till $T = 15$. Fig.\ref{f-12} shows admissible set computed at different noise amplitude in this case. The comparison with simulations done on $(0,T)$ for $T = 40$ (see Fig.\ref{f-10}a) shows that with relatively small number of samples we have identified one of the two parameters ($\mathrm{Da}_d$ in this case) with a reasonable accuracy. We can use this as a starting point for a next stage, considering smaller set $G$, (for which  
$0.03 \leq \mathrm{Da}_d \leq 0.06$). 

\section{Conclusion} 
Deterministic and stochastic, single stage and multistage algorithms for identification of unknown adsorption and desorption rates in Langmuir isotherm are presented in conjunction with pore scale simulation of reactive flow. Breakthrough curves are the extra data needed for the solution of the parameter identification problem. Exact and noisy data are considered. In the latter case their impact on the identification procedure is discussed.

\begin{enumerate}
 \item 
The 2D mathematical model of the direct problem includes steady state Stokes equations, and convection--diffusion equation supplemented with Robin type boundary conditions accounting for adsorption and desorption. Henry and Langmuir isotherms describe the kinetics. A simple pore scale geometry described by periodic arrangement of cylindrical obstacles, is considered for illustration of the identification procedure. The key dimensionless parameters are specified. The approach is applicable for wide range of microgeometries and process parameters.

 \item 
The numerical solution is based on triangular grids and FEM with Taylor and Hood elements. Computations on series of refined grids are performed to confirm the convergence.

 \item 
Mass transport is simulated for as given velocity field (one way coupling). The numerical solution is based on FEM with piecewise linear elements. Crank-Nikolson scheme is used in the time discretization. Sensitivity studies are carried out to investigate the influence of different parameters on the reactive transport through the porous media.

 \item 
Deterministic and stochastic, single stage and multistage identification procedures are presented and tested. The influence of the noise in the data on the accuracy of the identification procedure is discussed. It should be pointed that in the real problem one can not identify single values for the seek parameters, instead, final admissible sets of parameters are identified. 
\end{enumerate} 

\section*{Acknowledgements}

VVG acknowledges the support by Mega-grant of the Russian Federation Government, N 14.Y26.31.0013. PNV acknowledges the support by Mega-grant of the Russian Federation Government, N 14.Y26.31.0013 and by the RFBR, N 17-01-00689.


\end{document}